\documentclass[12pt]{article}
\usepackage{amsmath,latexsym,amssymb}
\usepackage{amsfonts,mathtools}

\newcommand{\re}{{\bf  R}}

\newcommand{\ree}{{\bf  R}}

\newcommand{\ddd}{{\rm d}}
\newcommand{\vecu}{{\bf u}}
\newcommand{\vecom}{{\bf \omega}}
\newcommand{\vecer}{{\bf e_r}}
\newcommand{\ptl}{{\partial}}
\newcommand{\om}{{\omega}}

\newcommand{\qed}{\mbox{}\hfill\raisebox{-2pt}{\rule{5.6pt}
{8pt}\rule{4pt}{0pt}}\medskip\par}

\newcommand{\be}{\begin{equation}}
\newcommand{\ee}{\end{equation}}
\newcommand{\beary}{\begin{eqnarray}}
\newcommand{\eeary}{\end{eqnarray}}
\newcommand{\ela}[1]{\label{eq:#1}}

\newcommand{\er}[1]{$(\ref{eq:#1})$}
 \newtheorem{theorem}{Theorem}

\begin{document}
\title{Exact solutions of fluid equations on a sphere}

\author{Sun-Chul Kim\thanks{ Department of Mathematics,
Chung-Ang University,  221 Heukseok-dong, Dongjak-ku, Seoul
156-756 Korea. (Corresponding author: kimsc@cau.ac.kr)} \and Habin Yim\thanks{ Department of Mathematics,
Chung-Ang University,  221 Heukseok-dong, Dongjak-ku, Seoul
156-756 Korea.} }

\maketitle

\begin{abstract}
Exact solutions of both the Navier-Stokes and Euler equations are found on the surface of a sphere. Under the assumption of a vanishing convection term, the flow of two oppositely rotating point vortices at the poles turns out to be the unique common solution.
\end{abstract}

\\


\section{Introduction}

Fluid flows on the surface of a sphere are directly applicable to atmospheric or ocean circulations on a macroscopic scale where the curvature effect is not negligible. Moreover, due to the compactness of the sphere, the total flow displays different properties compared to the planar case. Many efforts have been made to understand and propose possible flow models, taking into account these aspects, which differ from two-dimensional planar fluid flow.

One such attempt is to find exact solutions that provide an accurate mathematical model of real physical flows and can also be used to improve the accuracy and reliability of numerical simulations. The main difficulty here is the high nonlinearity of fluid equations (Navier-Stokes, Euler), which permits no closed-form general solutions. However, various types of exact solutions have been found for two- or three-dimensional cases \cite{Berker, DrRi, Wang, AlCo, Kim24}, with or without certain symmetries.

For flows on the surface of a sphere, the situation is quite similar. The geometric and topological characteristics of a sphere, such as positive curvature and compactness, as well as nonlinearity, make it more difficult to find exact solutions in general. To the best of the author's knowledge, there are only a few known exact solutions for flows on the surface of a sphere \cite{Maetal, CrCl, CoKr}.

In connection with this, the first author specified and classified certain exact solutions of both the Navier-Stokes and Euler equations in two dimensions, as well as for the axisymmetric case \cite{KimOk06}. These are characterized by the vanishing nonlinear convection term in the equation and are possibly called {\em basic solutions}, as they include well-known fundamental solutions, such as Poiseuille and Couette flows, rigid rotation, etc., in fluid dynamics.

In this note, using a similar method, we attempt to find and classify analogous basic solutions of fluid flows on the surface of a sphere. More specifically, we consider two-dimensional stationary solutions of the Navier-Stokes (and Euler) equations on a sphere of radius
$R$. In coordinate free form, the equations are
\begin{eqnarray}
(\vecu\cdot\nabla)\vecu &=&  \nu \triangle \vecu - \frac{1}{\rho}\nabla p,
\label{eq:ns01}  \\
{\rm div}\, \vecu &=& 0.  \label{eq:ns02}
\end{eqnarray}
Here $\vecu $ denotes the velocity field, $p$ the pressure. $\nu$
and $ \rho$ are positive constants representing the kinematic
viscosity and the mass density, respectively.

If we confine the flow on the surface of sphere $\Sigma$, we define the
stream function $\psi = \psi (\theta, \phi)$ in the spherical
coordinate $(r, \theta, \phi),\, r >0, 0 <\theta < \pi, 0 < \phi \le 2\pi$  and represent the velocity $\vecu$
by,
$$ \vecu = \nabla \psi \times \vecer = ( 0, u_\theta, u_\phi) = \left( 0,
\frac{1}{R\sin \theta}\frac{\ptl \psi}{\ptl \phi},
-\frac{1}{R}\frac{\ptl \psi}{\ptl \theta} \right). $$ Taking curl to the
equation (\ref{eq:ns01}) we introduce the vorticity $ \vecom =
\nabla \times \vecu = \omega \vecer = ( \omega, 0, 0).$ For
convenience, we assume $R=1$ hereafter.

By these functions $\psi$ and $\om$, whose existence as
single-valued functions is assumed, the equations (\ref{eq:ns01})
and (\ref{eq:ns02}) are rewritten as follows:
\begin{eqnarray}
\frac{1}{R^2 \sin \theta}\left(\frac{\partial \psi}{\partial
\phi}\frac{\partial \omega}{\partial \theta} - \frac{\partial
\psi}{\partial \theta} \frac{\partial \omega}{\partial
\phi}\right)
 &=&  \nu \tilde{\nabla}^2
\omega, \ela{eq:ns03}  \\
- \tilde{\nabla}^2  \psi &=& \omega.  \label{eq:ns04}
\end{eqnarray}
where we define the Laplace-Beltrami operator $\tilde{\nabla}^2$
on the sphere by
\be
\tilde{\nabla}^2
\omega = \frac{1}{\sin\theta}\frac{\ptl}{\ptl \theta} \left( \sin
\theta \frac{\ptl \omega}{\ptl \theta} \right) + \frac{1}{\sin^2 \theta}
\frac{\ptl^2 \omega}{\ptl \phi^2}. \ela{lap_bel}
\ee

We are interested in the exact solutions and their physical
implications here. Before doing so, we list one important integral
constraint of the solution. As the sphere is a closed (and
compact) surface, the total vorticity should vanish by the Stokes'
theorem:
\be \int\!\!\int_{\Sigma} \omega \; dA = 0, \ela{eq:constr}
\ee
which is often called the ``Gauss constraint''.

\section{Solutions depending only on $\theta$}
We first seek basic solutions by assuming $\psi$ (and $\omega$)
is a function of $\theta$ only. From the equations \er{eq:ns03} and \er{lap_bel}, we derive
\[ \frac{d \omega}{d\theta} = \frac{K_1}{\sin\theta}, \] which
then is integrated into
\[ \omega(\theta)= K_1 \log \left( \tan\frac{\theta}{2} \right) + K_2, \] where $K_1, K_2$
are arbitrary constants. This solution has singularity at $\theta
= 0, \pi$ and thus undefined at the north and south poles. However,
for $K_2=0$, it satisfies the integral relation \er{eq:constr}
and the corresponding streamfunction and velocity are well defined
everywhere. Let us verify this fact. The streamfunction satisfies

\[ \frac{\ptl}{\ptl \theta} \left( \sin \theta \frac{\ptl
\psi}{\ptl  \theta} \right)  = -K_1 \log \left( \tan\frac{\theta}{2} \right), \] which
simplifies into
\[ \psi '' + (\cot \theta) \psi' = -K_1 \log \left( \tan\frac{\theta}{2} \right). \]
Solving for $u_\phi = -\psi'$ under the continuity on $u$ as
$\theta \rightarrow 0, \pi$,  we obtain
\[ u_\phi(\theta) = -\frac{K_1}{\sin \theta} \int_0^\theta \sin \varphi
\log \left( \tan\frac{\varphi}{2} \right) d\varphi. \] For convenience, putting $K_1=1$, the velocity
assumes its maximum value at
$\theta = \frac{\pi}{2}$ and monotonically decreases to zero at
the south and north poles. The vorticity then satisfies $$
\int_0^{2\pi} \!\!\int_0^{\frac{\pi}{2}} \omega \; d\theta d\phi =
- \int_0^{2\pi} \!\!\int_{\frac{\pi}{2}}^\pi \omega \; d\theta
d\phi = \log 2.
$$
This solution is symmetric with respect to the equator and is
independent with the Coriolis effect.
We call this the {\em basic} solution,
which identifies with two point vortices at the north and south poles with opposite vortex strengths.

\section{Solutions of $\psi= f(\omega)$ and $\omega=g(\psi)$}

Next, we seek exact solutions which have a functional relation
between $\psi$ and $\omega$ and vanishing convection term. We first suppose that the vector
field given by a stream function $\psi$ represents an Euler flow
and, at the same time, a Navier-Stokes flow. This amounts to
assuming that
\begin{equation}
\frac{\partial \psi}{\partial \phi}\frac{\partial
\omega}{\partial \theta} - \frac{\partial \psi}{\partial \theta}
\frac{\partial \omega}{\partial \phi} = 0 \label{eq:ns08}
\end{equation}
and
\begin{equation}
\tilde{\nabla}^2 \omega = 0. \label{eq:ns09}
\end{equation}
In the two dimensional flow case, these equations provide the most
fundamental solutions containing Poiseulle, Couette, rigid
rotation, and shear flows. Analogous examples are studied and listed in
the three dimensional axisymmetric case \cite{KimOk06}. In fact, it is a slightly over-determined system
and there appear some nontrivial solutions. However, in the case of a sphere surface, we state a negative result.

\begin{theorem}
There is no nontrivial solution except the basic solution satisfying the functional relation
$ \psi = F(\omega)$ with a smooth single-valued function $F$ on a
sphere.
\end{theorem}

{\em Proof.\;}
We first introduce a new variable $ \chi =
\log(\tan\frac{\theta}{2}) $ and regard the vorticity as a
function with $\chi$ and $\phi$ i.e. $\omega = \omega(\chi,
\phi).$  This simplifies the Laplace-Beltrami operator in the
usual Laplacian form : \[ \tilde{\nabla}^2 = \frac{1}{\sin^2
\theta} \left( \frac{\ptl^2}{\ptl \chi^2} + \frac{\ptl^2}{\ptl
\phi^2} \right)  = \frac{1}{\sin^2 \theta} \triangle_{\chi, \phi}.
\] We compute then by (\ref{eq:ns09})
$$
-\sin^2\theta \; \omega = \triangle_{\chi, \phi} \psi =
F^{\prime\prime}(\omega) | \nabla \omega|^2 + F^{\prime}(\omega)
\triangle_{\chi, \phi} \omega = F^{\prime\prime}(\omega)
 | \nabla \omega|^2.
$$
It  therefore follows that we may define $\Phi$ by
$$
\frac{1}{\sin^2 \theta} | \nabla \omega|^2  =
\frac{-\omega}{F^{\prime\prime}(\omega)} \equiv \Phi(\omega).
$$
Since $\omega$ is harmonic in $\chi$ and $\phi$, there exists a
complex analytic function $W(z) = W(\chi+{\rm i} \phi)$  such that
$ \omega = {\rm Re}\left[ W(z) \right]$. Consequently,  \[ \log
|\nabla \omega|^2  =\log\Phi + 2 \log (\sin\theta) = 2 {\rm
Re}\left[ \log \left( \frac{\ddd W}{\ddd z} \right) \right]\] is
harmonic, too. This implies that
$$
0 = \triangle_{\chi,\phi} (\log \Phi(\omega) + 2 \log
(\sin\theta)) = \left( \frac{\Phi^{\prime}}{\Phi} \right)^{\prime}
| \nabla \omega|^2 -2 \sin^2\theta,
$$
which gives us the differential equation for $\Phi$,
\begin{equation}
 \left( \frac{\Phi^{\prime}(\om(\theta,\phi)}{\Phi (\om(\theta,\phi)} \right)^{\prime}
  \Phi(\om(\theta,\phi) = 2.
\label{eq:gg}
\end{equation}
We will show that $\omega$ is in fact a function of $\theta$ only.
(Strictly speaking, $\omega = {\rm constant}$ in some open set and
(\ref{eq:gg}) holds in the remaining region. But, for the moment,
we proceed as if we have proved that $\omega$ is constant
throughout $\re^2$ or (\ref{eq:gg}) throughout $\re^2$.) The
equation (\ref{eq:gg}) can be easily integrated and we have $
\Phi(\omega) = Ae^{B\omega}$, where $A$ and $B$ are real
constants.  This, in turn, implies that
\begin{equation}
| \nabla \omega |^2 = \left| \frac{\ddd W}{\ddd z} \right|^2 /
\sin^2 \theta  = Ae^{B\omega }. \label{eq:kam04}
\end{equation}
Note that the case where $ \omega = {\rm constant}$ is included in this
equation as a special case that $A=0$.
Clearly $A \ge 0$.  Since $A=0$ implies the case where the
vorticity is constant, we may assume that $ A > 0$.  Applying
$\log$
 on (\ref{eq:kam04}),  we have
$$
 2 {\rm Re}\left[ \log \left(  \frac{\ddd W}{\ddd z} \right) \right]
= B  {\rm Re}\left[ W(z) \right] + \log A + 2 \log(\sin\theta).
$$
To find such solution, since $\sin \theta = \frac{2
e^\chi}{1+e^\chi}$, we should find an analytic function $f(\chi +
i \phi)$ satisfying the relation \[ \log\frac{2 e^\chi}{1+e^\chi}
= {\rm Re} f(\chi + i \phi), \] but there is no such function from the
Cauchy-Riemann condition.

We thus
conclude $\omega = \omega(\theta)$ and the solution then
degenerates to the basic solution case which we mentioned earlier.

\qed

The reverse case $\omega=g(\psi)$ is similar. The problem is to
find $\psi$ and $G$ such that
\begin{eqnarray}
- \triangle \psi &=& G(\psi),  \label{eq:goert01} \\
\triangle^2 \psi &=& 0.  \label{eq:goert02}
\end{eqnarray}

\begin{theorem}
If $ G : \ree \rightarrow \ree$ is a single-valued smooth
function, and
 {\rm  (\ref{eq:goert01})} and {\rm (\ref{eq:goert02}) }
are satisfied, then $\psi$ is again the basic solution.
\end{theorem}

{\em Proof.\;}
We have
$$
0 = \triangle G(\psi) = G^{\prime\prime}(\psi)|\nabla \psi|^2 +
G^{\prime}(\psi) \triangle \psi,
$$
which implies that
$$
|\nabla \psi|^2 =
\frac{G(\psi)G^{\prime}(\psi)}{G^{\prime\prime}(\psi)}.
$$
Since $ \nabla \omega = G^{\prime}(\psi)\nabla \psi$, we obtain
$$
|\nabla \omega|^2 = \frac{G(\psi)\left( G^{\prime}(\psi)\right)^3}
{G^{\prime\prime}(\psi)}.
$$
Let the right hand side be denoted by $ J(\psi)$.  Then, the same
procedure as in the previous case produces
$$
0 = \triangle( \log J(\psi)+ 2\log\sin\theta) =
\frac{J^{\prime}}{J} \triangle \psi + \left( \frac{J^{\prime}}{J}
\right)^{\prime} | \nabla \psi|^2 -\sin^2\theta,
$$
which can be written as follows:
$$
- \frac{J^{\prime}}{J} G + \left( \frac{J^{\prime}}{J}
\right)^{\prime} \frac{GG^{\prime}}{G^{\prime\prime}} =
\sin^2\theta.
$$
Again, we proceed by a similar argument as in the case
$\psi=f(\om)$. Namely, as the left side is a function of only $\psi$, let us write it into
$$
F(\psi(\theta, \phi)) =
\sin^2\theta.
$$
Differentiating on $\phi$ gives,
$$
F'(\psi(\theta, \phi)) \frac{\ptl \psi}{\ptl \phi} =0,
$$
so either $F'=0$ or $\frac{\ptl \psi}{\ptl \phi} =0$.
The first does not happen and it follows that $\psi = \psi(\theta)$,
again returning to the basic case. We found no new solution
for both cases.

\qed

\section{Concluding remark}

There are various exact solutions of the Navier-Stokes and Euler equations. Among them, we may characterize and classify basic solutions such as Couette flow, Poiseuille flow, solid body rotation, etc., as those with a vanishing advection term. In other words, they are common solutions of the Navier-Stokes and Euler equations. Such solutions exist in both two and three dimensions \cite{KimOk06}. However, on the surface of a sphere, we found only one solution involving two point vortices, which has singularities at two antipodal points. Aside from this, we analytically showed that no additional solutions exist. This is a negative result for the existence of exact solutions, imposing a rather strict condition for such common solutions. We suppose that this is likely due to the strong restrictions imposed by compactness and a positive curvature.

\section*{Acknowledgments}
This research was supported by the Chung-Ang University Research grant in 2024 and by the National Research Foundation of Korea(NRF) grant funded by the Korea government(MSIT) (No. RS-2023-00240538).

\end{document}